\documentclass[11pt]{article}   
\usepackage{amsthm, amsmath, amssymb}   
\usepackage{epsfig}
\usepackage[all]{xy}
\theoremstyle{plain}
\newtheorem{teorema}{Theorem}[section]
\newtheorem{lema}[teorema]{Lemma}

\theoremstyle{definition}
\newtheorem{definicion}[teorema]{Definition}
\newtheorem{ejemplo}[teorema]{Example}

\theoremstyle{remark}
\newtheorem{nota}[teorema]{Remark}

\newenvironment{dem}{\noindent{\it Proof}.}{\qed}

\newcommand{\K}{\wh{K}}
\newcommand{\got}{\mathfrak}

\renewcommand{\v}{\wh{v}}

\newcommand{\wh}{\widehat}
\newcommand{\om}{\omega}

\newcommand{\al}{\alpha}

\newcommand{\Ga}{\Gamma}
\newcommand{\De}{\Delta}

\newcommand{\vp}{\varphi}

\newcommand{\m}{{\got m}_v}

\renewcommand{\d}{\Delta_v}

\newcommand{\lcor}{{\rm [\kern - 1.8pt [}}
\newcommand{\rcor}{{\rm ]\kern - 1.8pt ]}}
\renewcommand{\[}{\left[}
\renewcommand{\]}{\right]}
\newcommand{\llcor}{\[\kern -2.5pt\[}
\newcommand{\rrcor}{\]\kern -2.5pt\]}

\newcommand{\ZZ}{{\mathbb{Z}}}

\newcommand{\CC}{{\mathbb{C}}}

\newcommand{\Y}{{\bf Y}}
\newcommand{\Z}{{\bf Z}}
\newcommand{\rango}{\mathrm{rank}}

\def\[[#1]]{[\![#1]\!]}
\def\olm{\leq_{\mathrm{lex}}}

\def\lex{\mathrm{lex}}
\def\lm{<_\lex}
\def\lM{>_\lex}
\let\em=\sl
\newcommand{\X}{{\bf X}}

\title{Monomial discrete valuations in $k\lcor\X\rcor$}
\author{Francisco J. Herrera Govantes \and Miguel \'Angel Olalla Acosta
  \and Jos\'e Luis Vicente C\'ordoba }
\date{July, 2003}

\begin{document}
\maketitle

\begin{abstract}
\cite{Bri2,Br-He} prove that all the rank one discrete valuations of
$k((X_1,X_2))$, centered in the ring $k\[[X_1,X_2]]$, come from the
usual order function, i.e. there exists a finite number of
transformations such that we obtain a new field $k((Y_1,Y_2))$ where
the lifting of $v$ is a monomial valuation given by $v(Y_1)=v(Y_2)=1$.

In this work we generalize this result to the rank $m$ discrete
valuation of $K=k((\X ))$, centered in $R=k\[[\X ]]$. We prove that, if
the dimension of $v$ is $n-m$, the maximum since \cite{Ab2}, then
there exists an inmediate extension $L$ of $K$ where the valuation is
monomial. Therefore we compute explicitly the residue field of the
valuation.
\end{abstract}

\section{Preliminaries.}

\begin{nota}
In this remark we remember some definitions and results given by
I. Kaplansky in \cite{Kap1}.

Let $K$ be a valued field, let $v$ be the valuation.
Let $L$ be an extension of $K$ and $v'$ a valuation that extends
$v$. We say that the extension $K\subseteq L$ is {\em immediate}
if the values group and the residue field of $v$ and $v'$ are the
same. We say that a valued field $K$ is {\em maximal} if it
doesn't admit proper immediate extensions.

A well ordered set $\{ a_i\}\subset K$, without last element in
$K$, is called {\em pseudo-convergent} if
$$v(a_j-a_i)<v(a_k-a_j)$$ for all $i<j<k$. Easily we have that, if
$\{ a_i\}$ is pseudo-convergent, then $v(a_j-a_i)=v(a_{i+1}-a_i),\
\forall i<j$. So we can use the abbreviation $\om_i$ for
$v(a_j-a_i)$, with $j>i$. Let us note that $\{\om_i\}$ is an
increasing set of elements of $\Ga$.

An element $a\in K$ is called {\em limit} of the pseudo-convergent
set $\{ a_i\}$ if $v(a-a_i)=\om_i$ for all $i$.

\cite{Kap1} prove the following results:

\begin{itemize}
\item If $K\subset L$ is an immediate extension, then every element
of $L\setminus K$ is a limit of any pseudo-convergent subset of
$K$ without limit in $K$.

\item A valued field $K$ is maximal if and only if there exists a
limit for all its pseudo-convergent subsets.
\end{itemize}

In \cite{Kru}, Krull shown the existence of, at least, one maximal
immediate extension for all valued field. Therefore, if $\K$ is
the completion of $K$ by the valuation $v$ and $\v$ is the unique
valuation of $\K$ that extends $v$, then $K\subset \K$ is an
immediate extension and $\K$ is a maximal valued field.

Kaplansky shown the unicity of the maximal immediate extension of
a valued field $K$ that satisfy a certain condition, that he calls
hypothesis A. In the case of zero characteristic this hypothesis A
is empty.
%
%
%


Finally, \cite{Kap1} gives an useful structure theorem for all
maximal valued field. If $\De$ is a field and $\Ga$ an ordered
abelian group, the set of all formal series $$\sum
a_it^{\al_i}\mbox{ con }a_i\in\De ,\ \al_i\in\Ga\mbox{ y
}\{\al_i\}\mbox{ well ordered},$$ is a field with the usual sum
and times. We shall denote this field by $\De (t^{\Ga})$. In $\De
(t^{\Ga})$ we can consider the valuation $\nu_t$ given by
$$\nu_t\left(\sum_{i\geq 1}a_it^{\al_i}\right) =\al_1\mbox{ with
}a_1\ne 0.$$ Krull shown in \cite{Kru} that, with this valuation,
$\De (t^{\Ga})$ is a maximal field.
\end{nota}

\begin{teorema}
(Kaplansky, 1942) Let $K$ be a maximal valued field with values
group $\Ga$ and residue field $\De$ that satisfies the hypothesis
A. Then $K$ is analytically isomorphic to the power series field
$\De (t^{\Ga})$.
\end{teorema}

\subsection{Notation and definitions.}

Let $K=k(({\bf X}))=k((X_1,\ldots ,X_n))$ be the quotient field of
the formal power series ring $R=k\[[\X]] =k\[[ X_1,\ldots ,X_n]]$.

\begin{nota}
\begin{itemize}
\item[1)] We shall write all series $f\in R$ as
$f=\sum_{A\in\ZZ_{0}^{n}}f_A {\bf X}^A$, where, if
$A=(a_1,\ldots,a_n)$, then ${\bf X}^A$ means $X_1^{a_{1}}\cdots
X_n^{a_{n}}$. We shall say $$ \mathcal{E}(f)=\{A\in\ZZ_{0}^n\mid
f_A\neq 0\} \, . $$

\item[2)] In $\ZZ^m$ we'll consider the lexicographic order,
it'll be denoted by $\olm$. This is a total order for the group
structure.

\item[3)] Let $0<m\leq n$ be an integer and
$$ L=\{B_1,\ldots,B_n\}\subset \ZZ_0^m\setminus\{0\} $$ such that
$L$ is a generator system of $\ZZ^m$. Each monomial ${\bf X}^A$ of
 $R$ has an element of $\ZZ_0^m$ associated, that is called its
$L$-{\em degree}\index{grado}, that is $$ \mathrm{degree}_{L}({\bf
X}^A)=\sum_{i=1}^n a_iB_i\, , \quad A=(a_1,\ldots,a_n)\, . $$
\end{itemize}
\end{nota}

\begin{definicion}{\label{12}} Let $0<m\leq n$ be an integer and
$$ L=\{B_1,\ldots,B_n\}\subset \ZZ_0^m\setminus\{0\} $$ such that
$L$ is a generator system of $\ZZ^m$. Let $v:\,
R\to\ZZ^m\cup\{\infty\}$ be the function such that $v(0)=\infty$
and $$ v(f)=\min_{\olm}\left\{\mathrm{degree}_{L}({\bf X}^A) \mid
A\in \mathcal{E}(f)\right\}\, , $$ with $f\neq 0$. The extension
of $v$ to $K|k$, which values group is $\ZZ^m$, is a rank $m$
discrete valuation, called monomial valuation associated to $L$.
\end{definicion}

Throughout this work, let $v$ be a rank $m$ discrete valuation of
$K|k$ centered in $R$. Let $R_v$, $\m$ and $\Ga =\ZZ^m$ be the
ring, the maximal ideal and the values group of the valuation $v$,
respectively. We'll denote by $\d =R_v /\m$ to the residue field
of $v$. Since \cite{Ab2} we know that the dimension of $v$ (the
transcendence degree of $k\subset \d$) is lesser or equal than
$n-m$. We shall suppose that the dimension of $v$ is the maximum,
$n-m$.

So in this work {\em valuation} means rank $m$ discrete valuation
of $K|k$ centered in $R$ and dimension $n-m$.

\begin{nota}
Let $\K$ be a maximal immediate extension of $K$. Since
\cite{Kap1} we can suppose that $\K$ is the completion of $K$ with
respect of $v$. Let $\v$ the only extension of $v$ to $\K$. We
know that there exists an analytic isomorphism of $\K$ in $\d
(t^{\Ga})$, so its restriction to the ring $R$ gives an injective
homomorphism $$\begin{array}{rcl} \vp :R=k\[[ {\bf X}]] & \to & \d
(t^{\Ga})\\ X_i & \mapsto & \sum_{j\geq 1}a_{i,j}t^{\al_{i,j}}
\end{array}$$
with $a_{i,j}\in\d$ y $\{\al_{i,j}\}\subset\Ga$ well ordered.

If we consider the extension of $\vp$ to the quotient field $K$
and the valuation $\nu_t$ of $\d (t^{\Ga})$ previously defined,
then $v=\nu_t\circ\vp$.
\end{nota}

The purpose of this work is to construct $\vp$ explicitly, in
order to obtain a parametric equation of $v$ and, in consequence,
a construction of the residue field of $v$, as an extension of the
field $k$.

Therefore we'll prove that for all valuation $v$ of $K|k$, there
exists an immediate extension $K\subset L=k((\Y ))$ such that the
valuation that extends $v$ is monomial.

In other words: Any valuation $v$ comes from a monomial valuation.
This result generalizes the obtained in \cite{Bri2,Br-He} for rank
one discrete valuations of $k((X_1,X_2))$.

\subsection{Monoidal Transformation and immediate extension.}

Let $v$ be a valuation of $K|k$. Let us consider the next monoidal
transformation in $K$: 
$$\begin{array}{rcl} 
k((\X )) & \to & L=k((\Y ))\\ 
X_i & \mapsto & Y_i\mbox{ if } i\ne 2\\ 
X_2 & \mapsto & Y_2Y_1
\end{array}$$
with $v(X_2)>_\lex v(X_1)$. Then we have the following theorem.

\begin{teorema}{\label{th16}}
With these conditions, the extension $K\subset L$ is immediate.
\end{teorema}

\begin{dem}
Let us consider the rings $R=k\[[\X ]]$ and $S=k\[[\Y ]]$, and the
diagram 
$$\xymatrix{ 
R \ar@{^{(}->}[rr]^-{\vp} \ar@{^{(}->}[d] & & \d (t^{\Ga})\\ 
S\ar[urr]_-{\psi}}$$ 
Where $\psi$ is the natural extension of $\vp$ to $L$,
i.e. $\psi(Y_2)=\vp(X_2)/\vp(X_1)$.

If $\psi$ is injective, then $v'=\nu_t\circ\psi$ is a valuation of
$L$ that extends $v$ and both has the same values group $\Ga$ and
the same residue field $\d$.

So we can suppose, by contradiction, that $\psi$ is not injective.
Let ${\got p}$ be the implicit ideal $\ker (\psi )$ and $L'$  the
quotient field of the ring $$\frac{k\[[\Y ]]}{\got p}.$$ Clearly
the restriction of $\psi$ to $L'$ is injective and its composition
with $\nu_t$ define a valuation, let us put $w_1$. So the
extension $K\subset L'$ is immediate.

Let $w_2$ be the valuation ${\got p}$--\'adic of $L$. This is a
discrete rank one valuation.

The composition of both valuations, let us put $w$, is a discrete
rank $m+1$ valuation of $L|k$ whose residue field is $\d$. Then we
have 
$$\rango (w)+\dim (w)=n+1>\dim k\[[\Y ]].$$ 
Since Abhyankar's theorem (\cite{Ab2}, Theorem 1, p. 330), we know 
$$\rango (w)+\dim (w)\leq\dim k\[[\Y ]],$$ 
so there is a contradiction.
\end{dem}

\begin{nota}
Here we are used the existence of such injective homomorphism, proved
by Kaplansky \cite{Kap1}. Later we'll give an explicit construction of
$\psi$. There isn't circular reasoning.
\end{nota}

The following example shows that the condition of maximal dimension of
$v$ is necessary.

\begin{ejemplo}
Let $R=\CC\[[ X_1,X_2,X_3]]$ and $K$ its quotient field. Let us
consider the injective homomorphism
$$\begin{array}{rcl}
\vp\colon R & \to & \CC (u)\[[ t]]\\
X_1 &\mapsto & t\\
X_2 &\mapsto & ut^2\\
X_3 &\mapsto & e^{ut}-1.
\end{array}$$
The composition of $\vp$ (really its extension to the quotient fields)
with the usual order function $\nu_t$ of $\CC (u)(( t))$ is a rank
one discrete valuation of $K$, named $v$. The dimension of $v$ is 1, the
transcendence degree of the extension $\CC\subset\CC (u)$.

If we make the monoidal transformation
$$\begin{array}{rcl} 
R=\CC\[[ X_1,X_2,X_3]] & \to & S=\CC\[[ Y_1,Y_2,Y_3]]\\ 
X_i & \mapsto & Y_i\mbox{ if } i=1,3\\ 
X_2 & \mapsto & Y_2Y_1
\end{array}$$
Then we have an homomorfism
$$\begin{array}{rcl}
\vp\colon S & \to & \CC (u)\[[ t]]\\
X_1 &\mapsto & t\\
X_2 &\mapsto & ut\\
X_3 &\mapsto & e^{ut}-1.
\end{array}$$
that is not injective.
\end{ejemplo}

\section{Constructing $\d$.}

\subsection{Basis of a subgroup of $\ZZ^m$.}

There are well known procedures to compute a basis of a subgroup
$\Ga_0\subset\ZZ^m$ knowing any set of generators. In this
subsection we describe an algorithm that will be very useful in
order to prepare the valuation.

Let $\{ A_1,\ldots ,A_n\}$ be the set of generators of
$\Ga_0\subset\ZZ^m$, with $A_i\lM {\mathbf 0}$ for all $i$. We can
suppose, without lost of generality, that $A_i\olm
A_j\ \forall i<j$. Let $A=(a_{i,j})\in {\cal M}_{n\times n}$ be the
matrix whose rows are the elements $A_i$.

We shall consider two transformations with the rows of the matrix $A$:

(1) $F_{i,j}(q)$: To change the row $i$ by itself plus $q$ times the
row $j$, with  $q\in\ZZ$.

(2) To interchange rows.

Clearly the group generated by the rows of the matrix $A$ is equal to
the group generated by the rows of any transformation of $A$.

%

\begin{nota} {\em Algorithm.}
The the matrix $A$ has the following echelon form

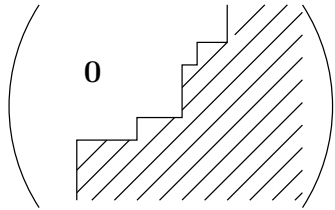
\begin{figure}[ht]
\begin{center}
\input{matriz.pstex_t}
\end{center}
\end{figure}

Let $j$ be the first column different of ${\mathbf 0}$ in $A$, let
$i$ be the first such that $a_{i,j}\ne 0$. Then we shall say that
$a_{i,j}$ is a {\em pivot}.

How $A_i\olm A_l\ \forall i<l$, clearly $a_{i,j}\leq a_{l,j}\
\forall i<l$. Let us put $a_{l,j}=q_la_{i,j}+r_l$, by making the
integer Euclidean division of $a_{l,j}$ by $a_{i,j}$.

We'll apply the following procedure with the first step, that is
perfectly exportable to the other steps:
\begin{itemize}
\item[I)] For each row $l$, with $i<l\leq n$, we make the following
  transformation:
\begin{itemize}
\item[a)] If $r_l\ne 0$ we do $F_{l,i}(-q_l)$. The new row $l$, that
  we denote by $A_l$ again, is such that $A_l\lm A_i$ and
  $0<a_{l,j}<a_{i,j}$.
\item[b)] If $r_l=0$ there are two possible situations:
\begin{itemize}
\item[1)] $A_l-q_lA_i\lM {\mathbf 0}$: In this case we make too the
  transformation $F_{l,i}(-q_l)$. The new row $l$ is such that $A_l\lm
  A_i$ and $a_{l,j}=0$. So, after a reordering $A_l$ raises a step.
\item[2)] $A_l-q_lA_i\olm {\mathbf 0}$: Then we make the change

  $F_{l,i}(1-q_l)$. The new row $l$, that we denote $A_l$ again, is
  such that $A_l\olm A_i$ and $m_{l,j}=m_{i,j}$. Let us remark that if
  $q_l=1$ then $A_l=A_i$, because we have suppose at begin that $A_i\olm A_l$.
\end{itemize}
\end{itemize}
\item[II)] After these transformations we reorder the rows. If every row of the
  first step of the matrix are equal, then we raise a step and we
  begin with the procedure.In other case, we apply again this algorithm.
\end{itemize}

Clearly we are doing the Euclidean algorithm in order to compute
the maximal common divisor of the elements $\{ a_{i,j},\ldots
,a_{n,j}\}$. Let $p_j$ be the maximal common divisor. By a finite
number of transformations we obtain a transformed matrix, that we
denote again by $A$, such that the pivot is equal to $p_j$ and the
last one divides all $a_{l,j}$ with $l\geq i$.

We have just arrived to situation I.b), so, by a finite number of
transformations, necessarily we must obtain a new matrix $A$ where
all the rows down the pivot $a_{i,j}=p_j$ are equals to $A_i$.

Hence, by applying this algorithm for each step of the matrix, we
obtain a matrix $B$ with only con $s$ different rows
$B_{i_1},\ldots ,B_{i_s}$ whose pivots are $p_{j_1},\ldots
,p_{j_s}$. Clearly $\Ga_0$ is isomorphic to
$p_{j_1}\ZZ\times\cdots\times p_{j_s}\ZZ$ and $\{ B_{i_1},\ldots
,B_{i_s}\}$ is a basis of $\Ga_0$.
\end{nota}

The algorithm described in this section allows us to prove the
following lemma, that we are going to use frequently for preparing
our valuation conveniently.

\begin{lema}{\label{lema}}
With the usual conditions over $K=k((\X ))$ and $v$, rank $m$
discrete valuation, if $\Ga_0$ is the subgroup generated by the
values of the elements $X_i$, the we can find, by a finite number
of monoidal transformations and interchanges of variables, an
immediate extension $L=k((\Y ))$ of $K$ such that each $Y_i$ has
the value in a basis $\{ B_1,\ldots ,B_s\}$ of $\Ga_0$.
\end{lema}

\begin{dem}
We have to apply the precedent algorithm to the matrix of the
values $v(X_i)$. Where $F_{l,i}(-q_l)$ means ``to apply $q_l$
monoidal transformations such that $X_l\mapsto Y_lY_i$'', and
reorder rows means ``to reorder variables according to its
values''.
\end{dem}

\subsection{Preparing $v$.}

As usually, we begin with a rank $m\leq n$ discrete valuation $v$ of
$K|k$, centered in the ring $R=k\[[\X ]]$. Let $\Ga =\ZZ^m$, $R_v$ and
$\m$ be the values group, the ring and the maximal ideal of the
valuation. $v$, respectively. We shall denote, as usual, by $\d$ to the
residue field of the valuation.

We'll suppose that the extension $k\subset\d$ is transcendent pure of
degree $\dim v=n-m$.


First, we apply the lemma \ref{lema} to obtain an immediate
extension $L=k((\Y ))$ of $K$ such that the set of values of the
elements $Y_i$ is a basis $\{ B_1,\ldots ,B_s\}$ of the subgroup
$\Ga_0$ generated by the values of the elements $X_i$.

{\em By convenience} we reorder the elements $Y_i$ in such way
that the first $s$ elements takes all the values of the basis
(i.e. $v(Y_i)=B_i$ for $i=1,\ldots , s$).

\begin{nota}
Let $A\in\Ga_0$ be such that $A=r_1B_1+\ldots
+r_sB_s$, then we shall denote by $R_A=(r_1,\ldots ,r_s,0,\ldots
,0)$ to the $n$--uple such that $v(\Y^{R_A})=A$.
\end{nota}

\subsection{The first transcendental residue.}

Our purpose is to give an explicit description of the injective
homomorphism $\psi\colon k\[[\Y ]]\to\d (t^{\Gamma})$, such that
$v=\nu_t\circ\psi$. In order not to complicate the exposition of this
construction, we are supposing that all the residues are in $k$ or
they are transcendental over the ground field. This condition seems
too strong, but after the proof of theorem \ref{th28} we'll explain
why this situation is really close to the general one.

\noindent
1) For the first elements we put $$\psi (Y_i)=t^{Bi}\ i=1,\ldots
,s.$$

\noindent
2) Let us take $Y_{s+1}$ the first element whose value is a linear
   combination of the elements $B_i$ (in fact it must be equal to some
   $B_i$). Let us suppose that $v(Y_{s+1})=B_{s+1,1}$, then we have
   two possibilities:
\begin{itemize}
\item[a)] For all $\al\in k$, $v(Y_{s+1}+\al\Y^{R_{B_{s+1,1}}})=
  B_{s+1,1}$. This fact means that the residue of
  $Y_{s+1}/\Y^{R_{B_{s+1,1}}}$ in $\d$ is transcendental over $k$. Let
  us put
  $$u_{s+1}=\frac{Y_{s+1}}{\Y^{R_{B_{s+1,1}}}}+\m$$
  and $\psi (Y_{s+1})= u_{s+1}t^{B_{s+1,1}}$.
\item[b)] There exists $\al\in k$ such that
  $v(Y_{s+1}+\al\Y^{R_{B_{s+1,1}}})\lM B_{s+1,1}$. So the residue of
  $Y_{s+1}/\Y^{R_{B_{s+1,1}}}$ in $\d$ is in $k$. Let us put
  $\al_{s+1,1}=\al$ and
  $v(Y_{s+1}+\al_{s+1,1}\Y^{R_{B_{s+1,1}}})=B_{s+1,2}>B_{s+1,1}$.
  There are two possibilities again:
\begin{itemize}
\item[i)] The new value $B_{s+1,2}\notin\Ga_0$. In this case, we make
  the change $Z_{s+1}=Y_{s+1}+\Y^{R_{B_{s+1,1}}}$, $Z_i=Y_i\ \forall
  i\ne s+1$ and go back to the beginning of the procedure by preparing the
  new valuation of $k((\Z ))$ with the lemma \ref{lema}.
\item[ii)] The value $B_{s+1,2}\in \Ga_0$. If there exists
  $\al_{s+1,2}$ such that
$$v\left(Y_{s+1}+\al_{s+1,1}\Y^{R_{B_{s+1,1}}}+
  \al_{s+1,2}\Y^{R_{B_{s+1,2}}}\right)=B_{s+1,3}\lM B_{s+1,2},$$
then we ask again if $B_{s+1,3}\in \Ga_0$. In the affirmative case,
if
  there exists $\al_{s+1,3}\in k$ such that
$$v\left(Y_{s+1}+\al_{s+1,1}\Y^{R_{B_{s+1,1}}}+
  \al_{s+1,2}\Y^{R_{B_{s+1,2}}}+
  \al_{s+1,3}\Y^{R_{B_{s+1,3}}}\right)=B_{s+1,4}\lM B_{s+1,3},$$
we go back to the beginning of the procedure.

We continue this procedure until we find a value
$B_{s+1,l}\notin\Ga_0$. If we cannot, there does not exist $\al_{s+1,l}\in k$
such that
$$v\left(Y_{s+1}+\sum_{k=1}^{l}\al_{s+1,k}\Y^{R_{B_{s+1,k}}}\right)=
B_{s+1,l+1}>B_{s+1,l}.$$ 
In the first case we make the change
$$Z_{s+1}=Y_{s+1}+\sum_{k=1}^{l-1}\al_{s+1,k}\Y^{R_{B_{s+1,k}}},\
Z_i=X_i\ \forall i\ne s+1$$ 
and move to the lemma \ref{lema}. In the second case we have that
the residue
$$u_{s+1}=\frac{Y_{s+1}+\sum_{k=1}^{l-1}
  \al_{s+1,k}\Y^{R_{B_{s+1,k}}}}{\Y^{R_{B_{s+1,l}}}}+\m$$
is transcendental over $k$. In this case we put
$$\psi (Y_{s+1})=\sum_{k=1}^{l-1}\al_{s+1,k}t^{B_{s+1,k}}
+u_{s+1}t^{B_{s+1,l}}.$$
\end{itemize}
\end{itemize}

We have to prove that this procedure ends up finding a new value that
it is not in $\Ga$ or a transcendental residue.

\begin{lema}
The situation described in 2.b.i) only occurs a finite number of times.
\end{lema}

\begin{dem}
There are two possible situations to find an element of value
$B\notin\Ga_0$:

1) The value $B$ is not a rational linear combination of the elements
of $\Ga_0$. In this case, the rank of the new subgroup of $\ZZ^m$
increases by 1. Trivially this fact only occurs a finite number of times.

2) The new value $B$ is a non integer rational linear combination of
the elements of the basis of $\Ga_0$. Let $\{ p_1,\ldots ,p_s\}$ be
the pivots that appears in the construction of the basis of
$\Ga_0$, let $\{ q_1,\ldots ,q_s\}$ be the ones of the new subgroup,
$\Ga_1$. As $\Ga_0\subset\Ga_1$, then $q_i\leq p_i\ \forall
i$ and, at least, one inequality is strict. As the {\em pivotes} are
greater or equal than 1, this only occurs a finite number of times.
\end{dem}

\begin{nota}{\label{nota}}
Let us suppose that we have a pseudo convergent set $\{ f_j\}$ of
elements of a valued field $K$, with value group $\ZZ^m$ with
lexicographic order. The set of values $\{\om_j\}$ of the elements
$f_j$ is a strictly increasing sequence of elements of the group $\ZZ^
m$. Let us suppose that the set $\{\om_j\}$ is bounded. Then, from
an index $j$ sufficiently big, we have
$$\om_j=(a_1,\ldots ,a_l,a_{l+1,j},\ldots ,a_{m,j}),$$ 
in such way that the first $l$ coordinates of the values $\om_j$ are
stabilized. Let us suppose that $l$ is the greatest integer between 1 and $m$
such that this fact occurs.

Let $f$ be a limit of  $\{ f_j\}$ and $\om =v(f)$ its value.
Then, if $\om =(b_1,\ldots ,b_n)$, as $\om\lM\om_j$ for all $j$,
then there exists an $l_0\leq l$ such that $b_i=a_i$ for all
$i=1,\ldots ,l_0-1$ and $b_{l_0}>a_{l_0}$.
\end{nota}

\begin{lema}{\label{lema2}}
The procedure described in the situation 2.b.ii) finds a value that
is not in  $\Ga_0$ or a transcendental residue over $k$.
\end{lema}

\begin{dem}
If the procedure is not finite, then we have a set of elements $\{
f_j\}_{j\geq 1}$ such that
$$f_j=Y_{s+1}+\sum_{l=1}^{j-1}\al_{s+1,l}\Y^{R_{B_{s+1,l}}}.$$
So $v(f_{j+1})=B_{s+1,j+1}\lM v(f_j)=B_{s+1,j}$ for all $j$.

Therefore, if $i<j<k$, we have
$$v(f_j-f_i)=B_{s+1,i}\lm v(f_k-f_j)=B_{s+1,j},$$
so $\{ f_j\}$ is a pseudo-convergent set.

Let us suppose that the set of values $\{ B_{s+1,j}\}$ is bounded.
Like in the remark \ref{nota}, let us suppose that, since an index
$j$ sufficiently great, the first $l$ coordinates of the
values are stabilized, where $l$ is the greatest integer between 1 and $m$
such that this fact occurs. Let us put $$v(f_j)=(a_1,\ldots
,a_l,a_{l+1,j},\ldots ,a_{m,j}).$$

Again by remark \ref{nota}, we know that any limit of $\{
f_j\}$ is such that its value is something like
$$(a_1,\ldots ,a_k,b_{k+1},\ldots ,b_m),$$ 
with $k<l$ and $b_{k+1}>a_{k+1}$. Let us take the series 
$$Y_{s+1}+\sum_{j=1}^{\infty }\al_{s+1,j}\Y^{R_{B_{s+1,j}}}.$$ 
As limit, by convenience we put 
$$g_1=\sum_{j=1}^{\infty }\al_{s+1,j}\Y^{R_{B_{s+1,j}}}$$ 
and $v(Y_{s+1}+g_1)=B_{s+1,1}^1$.

If $B_{s+1,1}^1\notin\Ga_0$, the we have finished. In other case,
if there exists $\al_{s+1,1}^1$ such that 
$$v\left( Y_{s+1}+g_1+\al_{s+1,1}^1\Y^{R_{B_{s+1,1}^1}}\right)
=B_{s+1,2}^1\lM B_{s+1,1}^1,$$ 
then we continue with our procedure. If it is infinite again and the
values contained are bounded, then we shall have a pseudo-convergent
set with limit
$$Y_{s+1}+g_1+\sum_{j=1}^{\infty }\al_{s+1,j}^1\Y^{R_{B_{s+1,j}^1}}.$$
Again by convenience we put
$$g_2=\sum_{j=1}^{\infty }\al_{s+1,j}^1\Y^{R_{B_{s+1,j}^1}}$$ and
$v(Y_{s+1}+g_1+g_2)=B_{s+1,1}^2\lM B_{s+1,1}^1$.

If this procedure does not find an element not in $\Ga_0$ or a
transcendental residue, we have a series
$$Y_{s+1}+g_1+g_2+g_3+\cdots $$ 
such that the set of the partial sums
$$\left\{ Y_{s+1}+\sum_{i=1}^jg_i\right\}_{j\geq 1}$$ 
is pseudo-convergent.

If the values of this set are bounded again, then, from one $j$
sufficiently great, the first $k$ coordinates of each value are
stabilized, where $k$ is the greatest integer between 1 and $m$ such that
this fact occurs. As all the values of this set are greater than
the first limit we have found, then it must be $k<l$.

Therefore, if we don't find a transcendental residue, after an
infinite procedure of calculation of limits of pseudo-convergent
sets, the values, if they are bounded, become stable in lesser and
lesser coordinates.

So this procedure must end finding an element that is not in
$\Ga_0$ or a new transcendental residue. The opposite fact is
equivalent to construct a series 
$$f=Y_{s+1}+g(Y_1,\ldots ,Y_s)\in k((\Y ))$$ 
such that the set of values of the partial sums is not a bounded
strictky increasing sequence. This means that $v(f)=\infty$. So $\psi
(f)=0$, but this is a contradiction because $\psi$ is injective.

So, at the end of this procedure, we shall have a series
$f=Y_{s+1}+g(Y_1,\ldots ,Y_s)\in k((\Y ))$ such that
$v(f)=B_{s+1}$ and, either $B_{s+1}\notin\Ga_0$, or there does not
exist $\al$ such that $v(f+\al\Y^{R_{B_{s+1}}})\lM B_{s+1}$. In
this last case we take $$u_{s+1}=\frac{f}{\Y^{R_{B_{s+1}}}}+\m$$
and put $$\psi (Y_{s+1})=g(t^{B_1},\ldots
,t^{B_s})+u_{s+1}t^{B_{s+1}}.$$
\end{dem}

\subsection{Construction of $\d$.}
We can suppose that we have an immediate extension $L=k((\Y ))$ of
$K$ such that:

1) $\{ B_1,\ldots ,B_m\}$ is a basis of  $\Ga =\ZZ^m$. This means
that the situation 2.b.i) is not going to appear any more and all
the values that we find are integer linear combination of the
basis.

2) We have applied the described procedure with the first $k>m$
elements $Y_i$ obtaining: $$\begin{array}{l} \psi (Y_i)=t^{B_i},\
i=1,\ldots ,m\\ \psi (Y_j)=\sum_{l\geq
1}\al_{j,l}t^{B_{j,l}}+u_jt^{B_j},\ j=m+1,\ldots ,k
\end{array}$$
such that $B_j\lM B_{j,l}\lM B_{j,l-1}$ y
$\{\al_{j,l}\}\subset\De_{j-1}=k(u_{m+1},\ldots
,u_{j-1})\subset\d$.

Our purpose is to describe the procedure with the variable
$Y_{k+1}$, that is analogous to the described previously to find
the first transcendental residue.

Let $v(Y_{k+1})=B_{k+1,1}$, there exists $\al_{k+1,1}\in\De_{k}$
such that $v(Y_{k+1}+\al_{k+1,1}\Y^{R_{B_{k+1,1}}})=B_{k+1,2}\lM
B_{k+1,1}$?

1) If the answer is affirmative, then we pore the question for
the new element $Y_{k+1}+\al_{k+1,1}\Y^{R_{B_{k+1,1}}}$.

2) If it is negative then we put $$\psi
(Y_{k+1})=u_{k+1}t^{B_{k+1,1}}$$ and go to the following variable.

\begin{lema}
This procedure ends finding a new transcendental residue, eventually after an infinite number of calculations.
\end{lema}

\begin{dem}
The reasoning is the same that in lemma \ref{lema2}, if there is
not such transcendental residue, then we can construct a
pseudo-convergent set that has a limit $$f=Y_{k+1}+g(Y_1,\ldots
,Y_m)\in k((\Y )),$$ of value $\infty$. This implies that $\psi
(f)=0$, and this is a contradiction
\end{dem}

With this procedure we prove the following theorem:

\begin{teorema}{\label{th28}}
The residue field of $v$ is $k(u_{m+1},\ldots ,u_n)$.
\end{teorema}

\begin{dem}
The constructions described in this work permit us to find one
immediate extension $L=k((\Y ))$ of $K$ such that the valuation
that extends $v$ to the field $L$ is defined by the composition of
$$
\begin{array}{rcl}
\psi :k((\Y )) & \to & k(u_{m+1},\ldots ,u_n)(t^{\Ga})\\ Y_i &
\mapsto & t^{B_i},\ i=1,\ldots ,m\\ Y_j &\mapsto & \sum_{k\geq
1}\al_{j,k}t^{B_{j,k}}+u_jt^{B_j},\ j=m+1,\ldots ,n
\end{array}
$$ with the valuation $\nu_t$ of $k(u_{m+1},\ldots
,u_n)(t^{\Ga})$. Clearly the residue field of this valuation is
$k(u_{m+1},\ldots ,u_n)$.
\end{dem}

\begin{nota} {\em The general case: considering algebraic residues.}
If we find non-trivial algebraic residues, then the procedure does
not change, and we obtain the residue field
$$\d =k (\zeta_{m+1},\ldots ,\zeta_n),$$
where $\zeta_i$ is a collection, eventually infinite, of algebraic
residues $\{\alpha_{i,k}\}$ and one transcendental residue $u_i$.

All of these algebraic extensions are finite; in other case we can do
the same as the proof of theorem \ref{th16} to obtain a valuation with
dimension greater than $n-m$.

It could be another problem in the general case, because in fact we are
constructing a representant field of $\d$ as an intermediate extension
of $k\subset\sigma (\d )\subset K$, where $\sigma$ is a section of the
natural homomorphism $R_v\to\d$. As the rings of power series are
completes, we can apply Hensel's lemma to know that we can find an algebraic
representant $\alpha\in K$ for every algebraic residue $\alpha +\m$.
\end{nota}

\section{Monomial valuations.}

Althought the number of operations to calculate $\psi$ explicitly
with the previous procedure is infinite, the number of monoidal
transformations, coordinates changes and interchanges of variables
that transform $K$ in $L$ is finite, because these only appear
when we find a value $B\notin \Ga_0$ and this fact occurs a finite
number of times. Following the trace of these transformations we
obtain a map $\vp :K\to\d (t^{\Ga})$ that is the restriction of $\psi$
and parametrizes $v$. In a similar way we can obtain some representatives
of every residue $u_i$ depending on $\X$.

Finally we give a theorem that generalizes the results obtained in
\cite{Bri2,Br-He} for rank one discrete valuations of $k((X_1,X_2))$
centered in $k\[[ X_1,X_2]]$.

\begin{teorema}
For all rank $m$ discrete valuation $v$ of $K|k$, with dimension $n-m$
and centered in $R$, there exists an immediate extension $L=k((\Z ))$
of $K$ such that the lifting of $v$ to $L$ is monomial.
\end{teorema}

\begin{dem}
In fact, by means of the previous procedure we obtain an immediate
extension $k((\Y ))$ of $K$ such that the lifting of $v$ is
$v=\nu_t\circ\psi$, with 
$$
\begin{array}{rcl}
\psi :k((\Y )) & \to & k(u_{m+1},\ldots ,u_n)(t^{\Ga})\\ Y_i &
\mapsto & t^{B_i},\ i=1,\ldots ,m\\ Y_j &\mapsto & \sum_{k\geq
1}\al_{j,k}t^{B_{j,k}}+u_jt^{B_j},\ j=m+1,\ldots ,n.
\end{array}
$$ 
This valuation is not monomial, but if we make the change
$$
\begin{array}{rcl}
k((\Y )) & \to & k((\Z ))\\ Y_i & \mapsto & Z_i,\ i=1,\ldots ,m\\
Y_j &\mapsto & Z_j-\sum_{k\geq 1}\al_{j,k}\Z^{R_{B_{j,k}}},\
j=m+1,\ldots ,n,
\end{array}
$$ 
then the natural extension of $\psi$ to $k((\Z ))$ is
$$
\begin{array}{rcl}
\phi :k((\Z )) & \to & k(u_{m+1},\ldots ,u_n)(t^{\Ga})\\ Y_i &
\mapsto & t^{B_i},\ i=1,\ldots ,m\\ Y_j &\mapsto & u_jt^{B_j},\
j=m+1,\ldots ,n
\end{array}
$$ 
and the valuation $v=\nu_t\circ\phi$ is monomial.
\end{dem}

\begin{nota}
In the general case, considering algebraic residues, we obtain an
inmediate extension $L=k'((\Z ))$, where $k'$ is an extension of $k$
contained in $\sigma (\d )$. The lifting of $v$ in $L$ is a monomial valuation.
\end{nota}


\begin{ejemplo}
Let $R=k\[[X_1,X_2,X_3,X_4 ]]$ be the power series ring, with $k=\ZZ
/\ZZ_5$. Let $K$ be the quotient field. Let us consider the series
$W_1=X_4-\sum_{i\geq 1}X_3^{3i}$ and $W_2=X_2-\sum_{i\geq
  1}X_1^i$. Let us consider the following valuations:

$v_1$: $(W_1)$--adic valuation of $K$,

$v_2$: $(W_2)$--adic valuation of $K/(W_1)$ and

$v_3$: monomial valuation of $K/(W_1,W_2)$ defined by the values
$v_3(X_1)=v_3(X_3)=1$.

Finally, let $v$ be the composite valuation $v=v_3\circ v_2\circ
v_1$. It is a rank 3 discrete valuation of $K|k$, such that: 
$$
\begin{array}{l}
  v(X_1)=v(X_3)=(0,0,1) \\
  v(W_2)=(0,1,0) \\
  v(W_1)=(1,0,0),
\end{array}$$
wherefrom
$$\begin{array}{l}
  v(X_2)=v(W_2+\sum_{i\geq 1}X_1^i)= (0,0,1) \\
  v(X_4)=v(W_1+\sum_{i\geq 1}X_3^{3i})= (0,0,3).
\end{array}
$$ 
So the values of the variables generate the subgroup $\Ga_0=\{
0\}\times\{ 0\}\times\ZZ$. From lemma \ref{lema} we know that the
transformation $X_4\mapsto Y_4Y_1^2$ takes $K$ in to an immediate
extension $L=k((Y_1,Y_2,Y_3,Y_4))$, such that $v(Y_i)=(0,0,1),\
i=1,2,3,4$. We shall apply the given procedure with this example to
obtain a parametrization, $\psi$, of $v$.

Let us put $\psi(Y_1)=t^{(0,0,1)}$. We shall begin with the
variable $Y_2$, as $v(Y_2)=(0,0,1)$, and we ask if there exists $\al\in
k$ such that $v(Y_2+\al Y_1)\lM (0,0,1)$. By the construction of $v$,
we know that $v(Y_2-Y_1)=v(X_2-X_1)= (0,0,2)\lM (0,0,1)$.
Following the procedure, we find that
$$v\left( Y_2-\sum_{i=
1}^liY_1^i\right) =(0,0,l+1)\lM (0,0,l).$$
So we have a pseudo-convergent set $\{
f_l\}$, with $f_l=Y_2-\sum_{i= 1}^liY_1^i$. A limit of this set is 
$$f_{\infty}=Y_2-\sum_{i\geq 1}iY_1^i.$$
As
$v(f_{\infty})=v(W_2)=(0,1,0)\notin\Ga_0$, we put
$$\psi (Y_2)=\sum_{i\geq 1}it^{(0,0,i)}+t^{(0,1,0)}.$$

We continue with $Y_3$, now $\Ga_0=\{
0\}\times\ZZ\times\ZZ$. Since it does not exist $\al\in k$ such that
$v(Y_3+\al Y_1)\lM (0,0,1)$, we put $$\psi (Y_3)=u_3t^{(0,0,1)}.$$

Let $\De_3=k(u_3)$, with $u_3=Y_3/Y_1+\m$. Let us go to the variable
$Y_4=X_4/X_1^2$, $v(Y_4)=(0,0,1)$.By the construction of $v$, we know
that
$$v\left( Y_4-u_3^3Y_1\right) =v\left(\frac{X_4-X_3^3}{X_1^2}\right)
=(0,0,4)\lM (0,0,1).$$
Following the procedure we have
$$v\left( Y_4-\sum_{i=
1}^lu_3^{3i}Y_1^{3i-2}\right) =v\left(\frac{X_4-\sum_{i=
1}^lX_3^{3i}}{X_1^2}\right) =(0,0,3i+1).$$
Wherefrom we come to have a pseudo-convergent set $\{ g_l\}$, with
$g_l=Y_4-\sum_{i=1}^lu_3^{3i}Y_1^{3i-2}$. A limit of this set is
$$g_{\infty
}=Y_4-\sum_{i\geq 1}u_3^{3i}Y_1^{3i-2}$$
and $v(g_{\infty })=v(W_1/X_1^2)=(1,0,-2)\notin\Ga_0$. We put
$$\psi
(Y_4)=\sum_{i\geq 1}u_3^{3i}t^{(0,0,3i-2)}+t^{(1,0,-2)}.$$

Therefore we have $\d =k((X_3/X_1)+\m )$ and
$$\left\{\begin{array}{l}
  \psi (X_1)=t^{(0,0,1)} \\
  \psi (X_2)= \sum_{i\geq 1}it^{(0,0,i)}+t^{(0,1,0)} \\
  \psi (X_3)=u_3t^{(0,0,1)} \\
  \psi (X_4)=\sum_{i\geq 1}u_3^{3i}t^{(0,0,3i)}+t^{(1,0,0)}.
\end{array}\right.
$$

This way, making the substitution $X_1=Z_1$, $X_2=Z_2+\sum
iZ_1^i$, $X_3=Z_3$ y $X_4=Z_4+ \sum Z_3^{3i}$, we have an immediate
extension $M=k((Z_1,Z_2,Z_3,Z_4))$ of $K$. The valuation that extends
$v$ to $M$ is the monomial valuation defined by
$v(Z_1)=v(Z_3)=(0,0,1)$, $v(Z_2)=(0,1,0)$ and $v(Z_4)=(1,0,0)$

\end{ejemplo}




\bibliographystyle{amsplain}
\bibliography{refer}

\providecommand{\bysame}{\leavevmode\hbox to3em{\hrulefill}\thinspace}
\providecommand{\MR}{\relax\ifhmode\unskip\space\fi MR }
\providecommand{\MRhref}[2]{%
  \href{http://www.ams.org/mathscinet-getitem?mr=#1}{#2}
}
\providecommand{\href}[2]{#2}
\begin{thebibliography}{1}

\bibitem{Ab2}
S.S. Abhyankar, \emph{On the valuations centered in a local domain.}, Amer. J.
  Math. \textbf{78} (1956), 321--348.

\bibitem{Bri2}
E.~Briales, \emph{Constructive theory of valuations.}, Commun. Algebra
  \textbf{17} (1989), no.~5, 1161--1177.

\bibitem{Br-He}
E.~Briales and F.J. Herrera, \emph{Construcci\'{o}n expl\'{\i}cita de las
  valoraciones de un anillo de series formales en dos variables.}, Actas X
  Jornadas Hispano-Lusas (Murcia), vol.~II, 1985, pp.~1--10.

\bibitem{Kap1}
I.~Kaplansky, \emph{Maximald fields with valuations}, Duke Math. J. \textbf{9}
  (1942), 303--321.

\bibitem{Kru}
W.~Krull, \emph{Algerneine bewertungstheorie}, J. Reine Angew. Math.
  \textbf{167} (1931), 160--196.

\end{thebibliography}

\end{document}